\documentclass[11pt]{article}

\usepackage{amsfonts}

\begin{document}

\def \a{{{\frak a}}}
\def \al{\alpha}
\def \ar{{\alpha_r}}
\def \A{{\mathbb A}}
\def \Ad{{\rm Ad}}
\def \b{{{\frak b}}}
\def \bs{\backslash}
\def \B{{\cal B}}
\def \cent{{\rm cent}}
\def \C{{\mathbb C}}
\def \CA{{\cal A}}
\def \CB{{\cal B}}
\def \CE{{\cal E}}
\def \CF{{\cal F}}
\def \CG{{\cal G}}
\def \CH{{\cal H}}
\def \CL{{\cal L}}
\def \CM{{\cal M}}
\def \CN{{\cal N}}
\def \CP{{\cal P}}
\def \CQ{{\cal Q}}
\def \CO{{\cal O}}
\def \det{{\rm det}}
\def \e{\epsilon}
\def \End{{\rm End}}
\def \Fx{{\frak x}}
\def \FX{{\frak X}}
\def \g{{{\frak g}}}
\def \ga{\gamma}
\def \Ga{\Gamma}
\def \h{{{\frak h}}}
\def \Hom{{\rm Hom}}
\def \Im{{\rm Im}}
\def \Ind{{\rm Ind}}
\def \k{{{\frak k}}}
\def \K{{\cal K}}
\def \la{\lambda}
\def \lap{\triangle}
\def \Lie{{\rm Lie}}
\def \m{{{\frak m}}}
\def \mathqed{\hfill \Box} 

\def \mod{{\rm mod}}
\def \n{{{\frak n}}}
\def \name{\bf}
\def \N{\mathbb N}
\def \ord{{\rm ord}}
\def \O{{\cal O}}
\def \p{{{\frak p}}}
\def \ph{\varphi}
\def \prf{{\bf Proof: }}
\def \qed{$ $\newline $\frac{}{}$\hfill {\rm Q.E.D.}\vspace{15pt}}
\def \Q{\mathbb Q}
\def \res{{\rm res}}
\def \R{{\mathbb R}}
\def \Re{{\rm Re \hspace{1pt}}}
\def \ra{\rightarrow}
\def \rank{{\rm rank}}
\def \Rep{{\rm Rep}}
\def \sign{{\rm sign}\hspace{2pt}}
\def \supp{{\rm supp}}
\def \t{{{\frak t}}}
\def \T{{\mathbb T}}
\def \tr{{\hspace{1pt}\rm tr\hspace{1pt}}}
\def \vol{{\rm vol}}
\def \V{{\cal V}}
\def \z{\frak z}
\def \Z{\mathbb Z}
\def \={\ =\ }

\newcommand{\rez}[1]{\frac{1}{#1}}
\newcommand{\der}[1]{\frac{\partial}{\partial #1}}
\newcommand{\binom}[2]{\left(  
\begin{array}{c}#1\\#2\end{array}\right)}
\newcommand{\norm}[1]{\parallel #1 \parallel}
\newcommand{\krn}[3]{\langle #1 | #2 | #3 \rangle}

\newcounter{lemma}
\newcounter{corollary}
\newcounter{proposition}
\newcounter{theorem}
\newcounter{zwisch}

\renewcommand{\subsection}{\refstepcounter{subsection}\stepcounter{lem 
ma} 

	\stepcounter{corollary} \stepcounter{proposition}
	\stepcounter{conjecture}\stepcounter{theorem}
	$ $ \newline
	{\bf \arabic{section}.\arabic{subsection}\hspace{8pt}}}

\newtheorem{conjecture}{\stepcounter{lemma} \stepcounter{corollary} 	
	\stepcounter{proposition}\stepcounter{theorem}
	\stepcounter{subsection}\hskip-12pt Conjecture}[section]
\newtheorem{lemma}{\stepcounter{conjecture}\stepcounter{corollary}	
	\stepcounter{proposition}\stepcounter{theorem}
	\stepcounter{subsection}\hskip-7pt Lemma}[section]
\newtheorem{corollary}{\stepcounter{conjecture}\stepcounter{lemma}
	\stepcounter{proposition}\stepcounter{theorem}
	\stepcounter{subsection}\hskip-7pt Corollary}[section]
\newtheorem{proposition}{\stepcounter{conjecture}\stepcounter{lemma}
	\stepcounter{corollary}\stepcounter{theorem}
	\stepcounter{subsection}\hskip-7pt Proposition}[section]
\newtheorem{theorem}{\stepcounter{conjecture} \stepcounter{lemma}
	\stepcounter{corollary}\stepcounter{proposition}		
	\stepcounter{subsection}\hskip-11pt Theorem}[section]

\title{Combinatorial $L^2$-determinants}
\author{{\small by}\\ {} \\ Anton Deitmar\\ {}\\ to appear in: Proc.  
Edinburgh Math. Soc.}
\date{}
\maketitle

\pagestyle{myheadings}
\markright{COMBINATORIAL $L^2$-DETERMINANTS}

\noindent
{\bf ABSTRACT:} We show that the zeta function of a regular graph  
admits a representation as a quotient of a determinant over a  
$L^2$-determinant of the combinatorial Laplacian.

\vspace{10pt}

\tableofcontents

\begin{center}
{\bf Introduction}
\end{center}

In \cite{det} we showed that a geometric zeta function of a locally  
symmetric space of rank one admits a representation as a quotient  
$\frac{\det(\lap +P(s))}{\det_\Ga(\lap +P(s))}$ where $P$ is a  
polynomial and $\lap$ a generalized Laplacian, det refers to the  
determinant and $\det_\Ga$ to the $L^2$-determinant.
The combinatorial counterpart to these geometric zeta functions are  
zeta functions of finite graphs.
In this paper we show an analogous formula for finite graphs.
For the sake of conceptual clarity we will give the full proof only  
for regular graphs, i.e. graphs of constant valency, but the methods  
are easily seen to cover the general case as well.
By the results of \cite{bass} we are reduced to a computation of the  
$L^2$-determinant of $\lap+\la$ which cannot be calculated directly  
because of the combinatorial complexity.
We circumvent this problem in using a technique developed in  
\cite{l2det} which essentially gives a way to compute  
$L^2$-determinants as limits of ordinary determinants.
See \cite{luck} for a similar assertion.
Using the results of \cite{bass} or \cite{hash} the proof of the main  
theorem (\ref{main}) then becomes easy.

If the valency of the graph is $q+1$ for $q$ being a prime power the  
same result can be derived by means of the harmonic analysis of the  
locally compact group $G=SL_2(F)$ for a $p$-adic field $F$, since  
then the graph may be compared to a quotient of the Bruhat-Tits  
building of $G$.
In \cite{padgeom} we generalized the notion of geometric zeta  
functions to higher rank Bruhat-Tits buildings.
The question whether there is an analogous formula in the higher rank  
case is up to now unanswered.

\section{The convergence Theorem}

Let $(Y,d)$ be a discrete metric space which is {\bf proper}, i.e.  
for any $y\in Y$ and any $R>0$ the set of all $x\in Y$ with  
$d(x,y)\le R$ is compact.
By the discreteness this means the latter set is finite.
Let $\Ga$ be a group of isometries of $Y$ which acts freely and such  
that the quotient $\Ga\bs Y$ is finite.
Note that the existence of $\Ga$ implies that $Y$ is {\bf uniformly  
proper}, i.e. for any $y\in Y$ and any $R>0$ the number of points $x$  
with $d(x,y)\le R$ is bounded by a constant only depending on $R$ but  
not on $y$.

Let $\CF$ be a sheaf of complex vector spaces over $Y$ of dimension  
$r$.
By the discreteness this just means that to each $y\in Y$ there is  
attached a $r$-dimensional vector space $\CF(y)$.
Assume the action of $\Ga$ lifts to $\CF$ such that it preserves  
stalks and is linear on each stalk.
For $R>0$ and $M\subset Y$ let $U_R(M)$ be the $R$-neighborhood of  
$M$, i.e. $U_R(M)$ is the set of $y\in Y$ with $d(y,M)<R$.

Let $\Ga(\CF)$ be the space of global sections of $\CF$ and let  
$\Ga_c(\CF)$ be the subspace of sections with compact support.
Let $R>0$, a linear operator $D:\Ga_c(\CF)\ra \Ga(\CF)$ is said to  
have {\bf finite propagation speed $\le R$} if for any  
$\ph\in\Ga_c(\CF)$ we have $\supp D\ph\subset U_R(\supp\ph)$.
An operator with finite propagation speed maps $\Ga_c(\CF)$ to  
itself.

\begin{lemma}
Suppose the operator $D$ maps $\Ga_c(\CF)$ to itself.
Suppose further that $D$ commutes with the action of $\Ga$ on  
$\Ga_c(\CF)$.
Then $D$ is of finite propagation speed.
\end{lemma}

\prf
Let $D$ be an operator stabilizing $\Ga_c(\CF)$.
A $\delta$-section at $y\in Y$ of the sheaf $\CF$ is a section  
$\delta$ which vanishes outside $\{ y\}$.
Let $F_D(y)$ be the maximum distance $d(x,y)$ of a point $x\in Y$  
with $D\delta(x)\ne 0$ for a $\delta$-section at $y$.
By the finite dimensionality of $\CF(y)$ this maximum is attained.
Note that $D$ has finite propagation speed if and only if the  
function $F_D$ is bounded.
Now assume $D$ to be $\Ga$-invariant, then so is $F_D$, which then,  
since $\Ga\bs Y$ is finite, must be bounded.
\qed

Since $\Ga$ acts freely on $Y$ the quotient $\CF_\Ga := \Ga\bs \CF$  
defines a sheaf on the finite set $Y_\Ga := \Ga\bs Y$.
The space of sections $\Ga(\CF_\Ga)$ can be identified with the space  
of $\Ga$-invariant sections $\Ga(\CF)^\Ga$ of $\CF$.
Note that the summation map $S:\Ga_c(\CF)\ra\Ga(\CF_\Ga)$ given by  
$S\ph(y) := \sum_{\ga\in \Ga}\ga^{-1}\ph(\ga y)$ is surjective.

\begin{lemma}
Let $D:\Ga_c(\CF) \ra \Ga(\CF)$ be linear, $\Ga$-invariant and of  
finite propagation speed, then there is an unique operator $D_\Ga$,  
the {\bf pushdown} of $D$ on $\Ga(\CF_\Ga)$ such that the diagram
$$
\begin{array}{ccc}
	\Ga_c(\CF)	&	{D\atop\ra}	&	\Ga_c(\CF)\\
	\downarrow S	&	{}		&	\downarrow  
S\\
	\Ga(\CF_\Ga) 	&	{D_\Ga\atop \ra}	& 	 
\Ga(\CF_\Ga)
\end{array}
$$
commutes.
\end{lemma}

\prf
Define $\tilde{D_\Ga}:\Ga_c(\CF)\ra \Ga(\CF_\Ga)$ by  
$\tilde{D_\Ga}(\ph) := SD(\ph)$.
The $\Ga$-invariance of $D$ and the finite propagation speed then  
implies that $\tilde{D_\Ga}$ vanishes on the kernel of $S$, hence  
induces an operator $D_\Ga$ as claimed.
\qed

Let $D:\Ga_c(\CF)\ra  \Ga(\CF)$ be linear, then $D$ has a kernel,  
i.e. it can be written as
$$
D\ph(x) \= \sum_{y\in Y} \langle x|D|y\rangle\ph(y),
$$where $\langle x|D|y\rangle\in\Hom_\C(\CF(y),\CF(x))$ and the sum  
is finite for each $x\in Y$.
It follows that $D$ has propagation speed $\le R$ if and only if  
$\langle x|D|y\rangle=0$ for all $x,y\in Y$ with $d(x,y)>R$.
Further $D$ is $\Ga$-invariant if an only if $\langle \ga x|D|\ga  
y\rangle =\ga \langle x|D|y\rangle\ga^{-1}$ for all $x,y\in Y$ and  
all $\ga\in \Ga$.

\begin{lemma}
Let $D$ be $\Ga$-invariant and of finite propagation speed then the  
operator $D_\Ga$ has the kernel:
$$
\langle x|D_\Ga |y\rangle \= \sum_{\ga\in \Ga} \ga^{-1}\langle \ga  
x|D|y\rangle.
$$
\end{lemma}

\prf
A computation.\qed

Suppose that the sheaf $\CF$ is hermitian and that $\Ga$ acts  
unitarily.
Let $L^2(\CF)$ be the space of square integrable sections of $\CF$.
Let $S$ be a set of representatives of $\Ga\bs Y$ then 

$L^2(\CF) \cong L^2(\Ga) \otimes L^2\left(\CF|_S\right) \cong  
L^2(\Ga)\otimes L^2(\CF_\Ga)$,
where $L^2(\Ga)$ is taken with respect to the counting measure and  
the tensor product is the tensor product in the category of Hilbert  
spaces.
Let $VN(\Ga)$ denote the von Neumann algebra of $\Ga$, i.e. the von  
Neumann algebra of all bounded operators on $L^2(\Ga)$ which commute  
with, say, the right action of $\Ga$.
It is easy to see that this algebra is topologically generated by the  
left translations $L_\ga$, $\ga\in \Ga$.
The algebra $VN(\Ga)$ carries a natural finite trace $\tau$ given by  
$\tau(\sum_{\ga\in \Ga}c_\ga L_\ga) = c_e$.
By the above we get that the von Neumann algebra $B(L^2(\CF))^\Ga$ of  
$\Ga$-invariant operators on $L^2(\CF)$ is isomorphic to
$B(L^2(\CF))^\Ga \cong VN(\Ga)\otimes B(L^2(\CF_\Ga))$.
Let $\tr_\Ga$ denote the trace on $B(L^2(\CF))^\Ga$ given by  
tensoring $\tau$ with the usual trace on $B(L^2(\CF_\Ga))$.

\begin{lemma}\label{trace}
Let $T\in B(L^2(\CF))^\Ga$ with kernel $\langle x|T|y\rangle$ then  
for any set $S$ of representatives of $\Ga \bs Y$ we have
$$
\tr_\Ga T \= \sum_{s\in S}\tr \langle x|T|x\rangle .
$$
\end{lemma}

\prf
As an element of $B(L^2(\CF))^\Ga\cong VN(\Ga)\otimes B(L^2(\CF|_S))$  
the operator $T$ writes as
$T\= \sum_{\ga\in \Ga} L_\ga \otimes PT|_S$,
where $P:\CF\ra\CF |_S$ is the projection.
Since the operator $PT|_S$ has kernel $\langle . |T| . \rangle  
|_{S\times S}$, the claim follows.
\qed

Let $\la >0$. For $\la$ large enough we may assume that the spectrum  
of the operator $D+\la$ lies in the right half plane.
Taking the standard branch of the logarithm the holomorphic  
functional calculus then allows to define the operator  
$(D+\la)^{-s}\in B(L^2(\CF))^\Ga$ for any $s\in \C$.
Set
$\zeta_{D+\la ,\Ga} (s) \ := \ \tr_\Ga((D+\la)^{-s})$,
and define the {\bf $L^2$-determinant} of $D+\la$ as
$$
\det_\Ga(D+\la) \ :=\ \exp(-\frac{d}{ds}|_{s=0} \zeta_{D+\la  
,\Ga}(s)).
$$

\begin{proposition}
The map $\la\mapsto \det_\Ga(D+\la)$ extends to a holomorphic  
function on the universal covering of $\C -{\rm Spec}(-D)$.
\end{proposition}

\prf
Let $X$ be the universal covering in question.
The map $\la\mapsto \log(D+\la)$ extends to an operator valued  
holomorphic map on $X$ and so does $(D+\la)^{-s} =  
\exp(-s\log(D+\la))$.
The claim follows.
\qed

A {\bf tower} of subgroups of $\Ga$ is a sequence $\Ga  
=\Ga_1\supset\Ga_2\supset\dots$ with $\cap_j\Ga_j =\{ 1\}$ and each  
$\Ga_j$ is normal of finite index in $\Ga$.

For any $\Ga$ and any natural number $N$ fix a standard $N$-th root  
of $\det(D_\Ga +\la)$ for $\Re(\la)>>0$ by requiring that
$\lim_{\la\ra+\infty}\frac{\det(D_\Ga +\la)^{1/N}}
			 {|\det(D_\Ga +\la)^{1/N}|}\= 1$.

We now come to the main result of this section.

\begin{theorem}\label{conv}
Let $D$ be $\Ga$-invariant and of finite propagation speed.
Let $(\Ga_j)$ be a tower in $\Ga$.
For $\Re(\la)>>0$ we have locally uniform convergence
$$
\det(D_{\Ga_j}+\la)^{\rez{[\Ga :\Ga_j]}}\ra \det_\Ga(D+\la)
$$
as $j\ra\infty$.
\end{theorem}

\prf
Let for $\Re(\la)>>0$ and $\Re(s)>>0$:
$$
F_j(\la ,s) \ := \ \frac{\tr(D_{\Ga_j}+\la)^{-s}}
		        {[\Ga :\Ga_j]}
	-\tr_\Ga (D+\la)^{-s}
$$ $$
\= \rez{\Ga(s)} \int_0^\infty t^{s-1} \left( \frac{\tr  
e^{-tD_{\Ga_j}}}
						  {[\Ga :\Ga_j]}
	-\tr_\Ga e^{-tD}\right) e^{-t\la} dt.
$$

\begin{lemma}\label{heat_kernel}
For any $z\in \C$ the operator $e^{zD_\Ga}$ has kernel
$$
\langle x |e^{zD_\Ga} |y\rangle \= \sum_{\ga\in \Ga} \ga^{-1}\krn{\ga  
x}{e^{zD}}{y}.
$$
\end{lemma}

\prf (Lemma)
The norms on $\CF(x)$ and $\CF(y)$ give rise to an norm on  
$\Hom_\C(\CF(y),\CF(x))$.
We have the estimate
\begin{eqnarray*}
\norm{\krn{x}{D_\Ga^n}{y}} &=& \sup_{\norm{v}=1}  
\norm{\krn{x}{D_\Ga^n}{y}v}\\
	&\le & \sup_{x,y,v} \norm{\krn{x}{D_\Ga^n}{y}v}\\
	&\le &\norm{D_\Ga^n}\le \norm{D_\Ga}^n,
\end{eqnarray*}
where $\norm{D_\Ga}$ is the operator norm on $L^2(\CF_\Ga)$.
It follows that we have absolute convergence in
\begin{eqnarray*}
e^{zD_\Ga}\ph(x) &=& \sum_{n\ge 0} \frac{z^n}{n!} D_\Ga^n\ph(x)\\
	&=& \sum_{n\ge 0}\frac{z^n}{n!} \sum_{y\in  
Y}\krn{x}{D_\Ga^n}{y}\ph(y)\\
	&=& \sum_{y\in Y} \sum_{n\ge 0}  
\frac{z^n}{n!}\krn{x}{D_\Ga^n}{y}\ph(y),
\end{eqnarray*}
so that
\begin{eqnarray*}
\krn{x}{e^{zD_\Ga}}{y} &=& \sum_{n\ge 0} \frac{z^n}{n!}  
\krn{x}{D_\Ga^n}{y}\\
	&=& \sum_{n\ge 0} \frac{z^n}{n!} \sum_{\ga\in\Ga}  
\ga^{-1}\krn{\ga x}{D^n}{y}\\
	&=& \sum_{\ga\in\Ga} \ga^{-1}\sum_{n\ge 0} \frac{z^n}{n!}  
\krn{\ga x}{D^n}{y}\\
	&=& \sum_{\ga\in\Ga}\ga^{-1}\krn{\ga x}{e^{zD}}{y}.
\end{eqnarray*}
The last equation follows from the above considerations which are  
clearly valid for trivial $\Ga$, too.
\qed

\begin{lemma}
There is a sequence $c_j>0$ tending to zero such that
$$
\left| \frac{\tr e^{zD_{\Ga_j}}}
	    {[\Ga :\Ga_j]}
	-\tr_\Ga e^{zD}\right| \le c_j|z|
$$
for all $z\in \C$ with $|z|\le 1$ and all $j\in \N$.
\end{lemma}

\prf
Let $S_j$ be a set of representatives of $\Ga_j\bs Y$.
By Lemma \ref{heat_kernel} we compute
\begin{eqnarray*}
\tr e^{zD_{\Ga_j}} &=& \sum_{s\in S_j} \tr  
\krn{s}{e^{zD_{\ga_j}}}{s}\\
	&=& \sum_{s\in S_j} \sum_{\ga\in\Ga_j }
		\tr(\ga^{-1}\krn{\ga s}{e^{zD_{\ga_j}}}{s}).
\end{eqnarray*}

Now let $S=S_1$ and assume that $S_j=\cup_{\sigma :\Ga /\Ga_j}\sigma  
S$, where $\sigma$ runs over a set of representatives of $\Ga  
/\Ga_j$.
The above equals
$$
\sum_{s\in S} \sum_{\sigma :\Ga /\Ga_j} \sum_{\ga\in\Ga_j}
		\tr(\ga^{-1}\krn{\ga \sigma s}{e^{zD_{\ga_j}}}{\sigma  
s}).
$$
Replacing $\ga$ by $\sigma\ga\sigma^{-1}$, which is possible since  
$\Ga_j$ is normal in $\Ga$, gives
\begin{eqnarray*}
&=& \sum_{s\in S} \sum_{\sigma :\Ga /\Ga_j} \sum_{\ga\in\Ga_j }
		\tr(\sigma \ga^{-1}\sigma^{-1}\krn{\sigma\ga   
s}{e^{zD_{\ga_j}}}{\sigma s})\\
&=& \sum_{s\in S} \sum_{\sigma :\Ga /\Ga_j} \sum_{\ga\in\Ga_j }
		\tr(\sigma \ga^{-1}\krn{\ga   
s}{e^{zD_{\ga_j}}}{s}\sigma^{-1}),
\end{eqnarray*}
since $e^{zD_{\ga_j}}$ is $\Ga$-invariant.
Now the $\sigma$-conjugation vanishes by taking traces, so we get
$$
 [\Ga :\Ga_j]\sum_{s\in S} \sum_{\ga\in\Ga_j}
	\tr(\ga^{-1}\krn{\ga  s}{e^{zD_{\ga_j}}}{s}).
$$

It follows that
$$
\frac{\tr e^{zD_{\Ga_j}}}
	    {[\Ga :\Ga_j]}
	-\tr_\Ga e^{zD}
\= \sum_{s\in S} \sum_{{\ga\in\Ga_j}\atop {\ga\ne 1}}
	\tr(\ga^{-1}\krn{\ga  s}{e^{zD_{\ga_j}}}{s}).
$$
In this sum we always have $\ga s\ne s$, so the diagonal of the  
kernel is never met.
Hence the term of $n=0$ in the sum $e^{zD} = \sum_{n\ge 0}  
\frac{z^n}{n!}D^n$ does not contribute.
It follows for $|z|\le 1$:
\begin{eqnarray*}
\left| \frac{\tr e^{zD_{\Ga_j}}}
	    {[\Ga :\Ga_j]}
	-\tr_\Ga e^{zD}\right|
&\le& \sum_{s\in S}\sum_{{\ga\in\Ga_j}\atop {\ga\ne 1}}
	\sum_{n\ge 1}\frac{|z|^n}{n!}
	|\tr(\ga^{-1}\krn{\ga  s}{D^n}{s})|\\
&\le& |z|\sum_{s\in S}\sum_{{\ga\in\Ga_j}\atop {\ga\ne 1}}
	\sum_{n\ge 1}\frac{1}{n!}
	|\tr(\ga^{-1}\krn{\ga  s}{D^n}{s})|,
\end{eqnarray*}
giving the claim.
\qed

To finish the proof of the theorem
we have to show that $F_j(\la ,s)$ tends to zero locally uniformly in  
$\la$ when $\Re(\la)>>0$ and $s$ in a neighborhood of zero.
To this end we split the integral into the pieces $\int_0^1$ and  
$\int_1^\infty$.
The second one converges for all $s$ when $\Re(\la)$ is large enough  
and the first one converges for $\Re(s)>-1$ by the last lemma.
Furthermore, this lemma shows the convergence in question for the  
part $\int_0^1$ already.

We now show that the second integral tends to zero.
For $\ph\in\Ga(\CF_{\Ga_j})$ let $\norm{\ph}_\infty := \sup_{y\in Y}  
\norm{\ph(y)}$ the sup-norm.
Let $\norm{D_{\Ga_j}}_\infty$ denote the corresponding operator norm.

\begin{lemma}
There is a constant $c>0$ such that $\norm{D_{\Ga_j}}_\infty\le c$  
for all $j\in\N$.
\end{lemma}

\prf
We estimate
\begin{eqnarray*}
\norm{D_{\Ga_j}\ph}_\infty &=& \sup_{x\in Y} \norm{D_{\Ga_j}\ph(x)}\\
&\le & \sup_{x\in Y}\sum_{y\in Y} \norm{\krn{x}{D}{y}\ph(y)}\\
&\le& \sup_{x\in Y}\sum_{y\in Y} \norm{\krn{x}{D}{y}}\norm{\ph(y)}.
\end{eqnarray*}

Now suppose $\norm{\ph}_\infty\le 1$ then we get
$$
\norm{D_{\Ga_j}\ph}\ \le \ \sup_{x\in Y}\sum_{y\in Y}  
\norm{\krn{x}{D}{y}}<\infty ,
$$
independent of $j$.
\qed

\begin{lemma}
There are $C_1,C_2 >0$s such that for all $j\in\N$:
$$
\left| \tr_\Ga e^{zD}\right| ,\left| \frac{\tr e^{zD_{\Ga_j}}}
	    {[\Ga :\Ga_j]}\right|
\ \le\ C_1e^{|z|C2}.
$$
\end{lemma}

\prf
We only show the assertion for $\left| \frac{\tr e^{zD_{\Ga_j}}}
	    {[\Ga :\Ga_j]}\right|$ since the other one is proven  
analogously.
Since for any operator $A$ on a finite dimensional euclidean space  
$V$ we have $|\tr A|\le (\dim V)\norm{A}_\infty$ it follows that 

\begin{eqnarray*}
|\tr e^{zD_{\Ga_j}}| &\le& (\# S)r[\Ga  
:\Ga_j]\norm{e^{zD_{\Ga_j}}}_\infty\\
&\le& (\# S)r[\Ga :\Ga_j]e^{\norm{{zD_{\Ga_j}}}_\infty},
\end{eqnarray*}
which by the last lemma implies the claim.
\qed

Now choose $\la$ with $\Re(\la)>C_2$ then it follows by the last  
lemma that the second integral of $F_j(\la ,s)$ converges dominatedly  
independent of $j$.
Hence it suffices that the integrand tends to zero pointwise.
This however is clear. The theorem is proven.
\qed

\section{Laplacians}

Let $X_\Ga =\Ga\bs X$ be a finite connected CW-complex with  
fundamental group $\Ga$ and universal covering $X$.
Let $\CF_\Ga$ be a hermitian locally constant sheaf of finite  
dimensional complex vector spaces and let $\CF$ be its pullback to  
$X$.
Choosing a basepoint $x_0$ gives a representation $\rho$ of $\Ga$ on  
the stalk over $x_0$.
For the cohomology we have $H^p(X_\Ga,\CF_\Ga) =H^p(\Ga ,\rho)$.
Let $q\ge 0$ and let $X_q$ be the set of $q$-dimensional cells of  
$X$.
We construct a discrete proper metric space $(Y,d)$ as follows:
The set $Y$ is given by $X_q$ and the distance $d(a,b)$ of two cells  
equals the minimal number of cells hit by a path joining a given  
point in $a$ to a given point in $b$ minus one.
For any cell $c$ of $X$ let $\CF(c)$ denote the stalk of $\CF$ at a  
given fixed point of $c$. (These points should be chosen  
$\Ga$-invariantly.)
This way $\CF$ induces a hermitian sheaf $\CF_Y$ on $Y$ on which  
$\Ga$ acts unitarily.
The combinatorial Laplacian $\lap_{q,\CF}$ now induces an operator on  
$\CF_Y$ which is $\Ga$-invariant and of finite propagation speed.
The theorem of the last section applies.

We will make use of this fact to prove a theorem on zeta functions of  
graphs.
So we restrict to the case $\dim X =1$, so $X$ is a tree and $X_\Ga$  
is a finite graph.
We give $X$ and $X_\Ga$ a metric by giving each edge the metric of  
the unit interval.
It then makes sense to speak of geodesics and especially of closed  
geodesics in $X_\Ga$.
Let $\ga$ be a closed geodesic and $x$ a point on its trace.
Since the sheaf $\CF$ is locally constant, parallel transport along  
$\ga$ gives a monodromy operator $m_\ga$ on the stalk $\CF(x)$.
We define the zeta function $Z_{\Ga ,\CF}$ of $\CF_\Ga$ as
$$
Z_{\Ga ,\CF}(u) \ :=\ \prod_{\ga}\det(1-u^{l(\ga)}m_\ga),
$$
where the product runs over all primitive closed geodesics, i.e.  
those, which are not a power of a shorter one.
The definition does not depend on the choices of points $x$ on the  
geodesics.
In general the product will be infinite but it is easy to show  
convergence for $u\in\C$ with $|u|$ sufficiently small.
In \cite{hash} it is shown that $Z_{\Ga ,\CF}$ extends to a rational  
function, indeed a polynomial, on $\C$.
Note that $\Ga$, being the fundamental group of a finite graph,  
contains a tower.

\begin{lemma}
Let $(\Ga_j)$ be a tower in $\Ga$ then $Z_{\Ga_j ,\CF}(u)^{\rez{[\Ga  
:\Ga_j]}}$ tends to $1$ for any $u\in \C$ with $|u|$ small.
\end{lemma}

\prf
Let $X^b_\Ga$ be the barycentric subdivision of $X_\Ga$, then  
$X^b_\Ga$ is bipartite, i.e. the vertex set of $X^b_\Ga$ can be  
written as a disjoint union $V=V_0\cup V_1$ where vertices in $V_i$  
are only connected by edges to vertices in $V_{1-i}$ for $i=0,1$.
Let $(Y,d)$ be the discrete proper metric space attached to the set  
of edges $X^b_1$ of $X^b$.
The sheaf $\CF$ can also be considered as a locally constant sheaf on  
$X^b$.
We will construct two operators $T_0$ and $T_1$ on $\CF_Y$.
Note first that for $a,b\in Y$, having a common vertex in $X^b$ then  
parallel transport gives an isomorphism $\ph_{a,b}:\CF(a)\ra\CF(b)$.
For $v\in\CF(a)$ let $\delta_v$ be the section of $\CF_Y$ that maps  
$a$ to $v$ and is zero elsewhere.
Then the $T_i$ are defined by
$$
T_i\delta_v := \sum_{{b\in Y}\atop{a\sim_i b}}\delta_{\ph_{a,b}(v)},
$$
where $a\sim_i b$ means that $a$ and $b$ have a common vertex in  
$V_i$ for $i=0,1$.
This prescription defines two $\Ga$-invariant operators of finite  
propagation speed on $\CF_Y$.
A computation shows
\begin{eqnarray*}
Z_{\Ga ,\CF}(u) &=& \det (1-u(T_0T_1)_\Ga)\\
&=& u^{2l}\det(\rez{u}-(T_0T_1)_\Ga),
\end{eqnarray*}
where $l$ is the number of edges of $X_\Ga$.
Theorem \ref{conv} implies
$$
Z_{\Ga_j ,\CF}(u)^{\rez{[\Ga :\Ga_j]}}\ra u^{2l}\det_\Ga\left(  
\rez{u}-(T_0T_1)_\Ga\right) ,
$$
as $j\ra \infty$ when $\Re(\rez{u})>>0$.
The operator $T_0T_1$ has the property that the diagonal of the  
kernel $\langle x|(T_0T_1)^n|x\rangle$ vanishes for any $n>0$.
This implies by Lemma \ref{trace} that $\det_\Ga(\rez{u}-(T_0T_1))  
=\det_\Ga(\rez{u})=u^{-2l}$, whence the claim.
\qed

We will say a graph is {\bf regular of valency $q+1$} for $q\in\N$,  
if every edge has two distinct endpoints and every vertex is  
connected to $q+1$ distinct edges.
To see which topological information is encoded in $Z_{\Ga ,\CF}$ we  
will connect it to the combinatorial Laplacian $\lap_{0,\CF_\Ga}$.

\begin{proposition}\label{hashi}
Let $\lap =\lap_{0,\CF_\Ga}$ and let $n$ be the number of vertices of  
$X_\Ga$.
Assume $X_\Ga$ is regular of valency $q+1$ for $q\ge 2$.
Then
$$
Z_{\Ga ,\CF}(u) = (1-u^2)^{-\chi(\rho)} u^n \det(\lap  
-(q+1-qu-\rez{u})),
$$
where $\chi(\rho)=\dim H^0(\Ga,\rho)-\dim H^1(\Ga,\rho)$ is the Euler  
characteristic.
\end{proposition}

\prf
In \cite{hash} it is shown that $Z_{\Ga ,\CF}(u)$ equals
$$
Z_{\Ga ,\CF}(u) = (1-u^2)^{-\chi(\rho)} \det(1-A_\rho u+qu^2),
$$
where $A_\rho$ is the adjacency operator of $\CF_\Ga$.
A calculation shows that $\lap =q+1-A_\rho$.
The claim follows.
\qed

The next theorem is the main result of this paper.

\begin{theorem}\label{main}
Let $\lap =\lap_{0,\CF_\Ga}$.
For $|u|$ small we have
$$
Z_{\Ga ,\CF}(u) \= \frac{\det(\lap +qu+\rez{u}-q-1)}
	  	        {\det_\Ga(\lap +qu+\rez{u}-q-1)}.
$$
\end{theorem}

\prf
At first note that, since the cohomology is computed by a finite  
dimensional complex it follows that
$$
\chi(\rho) \= \frac{(1-q)}{2}n\dim\rho,
$$
where $n$ is the number of edges on $X$.
Next, when $\Re(\rez{u})$ tends to infinity then so does $\Re(\la_u)$  
with $\la_u = qu+\rez{u}-q-1$.
Fix a tower $(\Ga_j)$.
For $\Re(\rez{u})>>0$ we have that $Z_{\Ga_j ,\CF}(u)^{\rez{[\Ga  
:\Ga_j]}}$ tends to $1$ as $j$ tends to infinity.
On the other hand, by Theorem \ref{conv} it follows that with $\lap_j  
=\lap_{0,\CF_{\Ga_j}}$ we have $\det(\lap_j +\la_u)^{\rez{[\Ga  
:\Ga_j]}}$ tends to $\det_\Ga(\lap +\la_u)$.
By Proposition \ref{hashi} we infer that 

$$
\det_\Ga(\lap +\la_u)\= (1-u^2)^{-\frac{q-1}{2}\dim(\rho)n}u^{-n}
$$
and thus the claim.
\qed

Using the results of \cite{bass} the theorem can be extended to  
arbitrary finite graphs, where the number $q$ has to be replaced by  
the operator $qf(x)=q(x)f(x)$, where $q(x)$ is the number of edges  
emanating from the vertex $x$.

\tiny
\hspace{-20pt}
Math. Inst.\\ INF 288\\ 69120 Heidelberg\\ GERMANY


\begin{thebibliography}{XXX}

\begin{scriptsize}

\bibitem{bass}
\bf Bass, H.:
\it The Ihara-Selberg zeta function of a tree lattice. 

\rm Int. J. Math. 3, No.6, 717-797 (1992).

\bibitem{det}
\bf Deitmar, A.:
\it A Determinant Formula for the generalized Selberg Zeta Function.
\rm Quarterly J. Math. 47, 435-453 (1996).

\bibitem{l2det}
\bf Deitmar, A.:
\it Regularized and $L^2$-Determinants.
\rm Proceedings of the London Mathematical Society 76, 150-174  
(1998).

\bibitem{padgeom}
\bf Deitmar, A.:
\it Geometric zeta functions on $p$-adic groups.
\rm Math. Japon. 47, No 1, 1-17 (1998).

\bibitem{hash}
\bf Hashimoto, K.:
\it Zeta functions of finite graphs and representations of p-adic  
groups.
\rm Automorphic forms and geometry of arithmetic varieties. Adv.  
Stud. Pure Math. 15, 211-280 (1989).

\bibitem{luck}
\bf L\"{u}ck, W.: 

\it Approximating $L^2$-invariants by their finite-dimensional  
analogues. 

\rm Geom. Func. Anal. 4 No. 4, 455-481, (1994).

\end{scriptsize}


\end{thebibliography}
\end{document}